\documentclass[11pt,reqno]{amsart}
\usepackage{amsmath,amsthm,amssymb,mathrsfs,stmaryrd}
\usepackage[all]{xy}
\usepackage{url}
\usepackage{geometry}
\usepackage[dvips]{color}
\usepackage{amsxtra}
\usepackage{amscd}
\usepackage{eucal}
\usepackage{epsfig}
\usepackage{graphics}

\usepackage{bm}

\usepackage[utf8]{inputenc}
\usepackage[OT1]{fontenc}

\usepackage{relsize}
\usepackage[bbgreekl]{mathbbol}
\usepackage{amsfonts}
\usepackage{textgreek}

\DeclareSymbolFontAlphabet{\mathbb}{AMSb}
\DeclareSymbolFontAlphabet{\mathbbl}{bbold}

\newcommand{\RR}{\mathbb{R}}

\newcommand{\cM}{\mathcal{M}}

\newcommand{\R}{\mathbb{R}}

\newcommand{\Q}{\mathbb{Q}}
\newcommand{\Z}{\mathbb{Z}}

\newcommand{\g}{\mathfrak{g}}

\newcommand{\nc}{\newcommand}
\nc{\on}{\operatorname}
\nc{\la}{\lambda}
\nc{\wh}{\widehat}
\nc{\wt}{\widetilde}
\nc{\sw}{{\mathfrak s}{\mathfrak l}}
\nc{\ghat}{\wh{\g}}
\nc{\hhat}{\wh{\h}}
\nc{\mc}{\mathcal}
\nc{\bi}{\bibitem}
\nc{\pa}{\partial}
\nc{\ppart}{(\!(t)\!)}
\nc{\pparl}{(\!(\la)\!)}
\nc{\zpart}{(\!(z^{-1})\!)}
\nc{\n}{{\mathfrak n}}
\nc{\ol}{\overline}
\nc{\mb}{\mathbf}
\nc{\bb}{{\mathfrak b}}
\nc{\su}{\wh\sw_2}
\nc{\h}{{\mathfrak h}}
\nc{\can}{\on{can}}
\nc{\ntil}{\wt{\n}}
\nc{\pone}{{\mathbb P}^1}
\nc{\bs}{\backslash}
\nc{\al}{\alpha}
\nc{\gt}{{\mathfrak g}'}
\nc{\ds}{\displaystyle}
\nc{\ep}{\varepsilon}
\nc{\alp}{\alpha}
\nc{\nab}{\nabla}
\nc{\RN}{\RR^n}
\nc{\RM}{\RR^m}
\nc{\Sc}{\mathrm{Schr}\ddot{\mathrm{o}}\mathrm{dinger}}
\nc{\Scd}{Schr\ddot{o}dinger}
\nc{\hb}{\hbar}
\nc{\st}{\section}
\nc{\sst}{\subsection}
\nc{\ta}{\theta}
\nc{\mm}{\mathrm}
\nc{\fr}{\frac}
\nc{\nb}{\nabla}
\nc{\pt}{\partial}
\nc{\ra}{\rightarrow}
\nc{\La}{\Leftarrow}
\nc{\Ra}{\Rightarrow}
\nc{\LRa}{\Leftrightarrow}
\nc{\lra}{\leftrightarrow}
\nc{\tim}{\times}
\nc{\oli}{\overline}
\nc{\uli}{\underline}
\nc{\ch}{\check}
\nc{\lt}{\left}
\nc{\rt}{\right}
\nc{\lap}{\Delta}
\nc{\infi}{\infty}
\nc{\intf}{\int_{-\infty}^\infty}
\nc{\de}{\delta}
\nc{\Gm}{\Gamma}
\nc{\gm}{\gamma}
\nc{\sse}{\subset}
\nc{\ssen}{\subsetneq}
\nc{\nsse}{\not\subset}
\nc{\lds}{\ldots}
\nc{\cds}{\cdots}
\nc{\cdt}{\cdot}
\nc{\dds}{\ddots}
\nc{\vds}{\vdots}
\nc{\tet}{\text}
\nc{\mr}{\mapsto}
\nc{\ml}{\mapsfrom}
\nc{\di}{\mm{dist}}
\nc{\til}{\tilde}
\nc{\ba}{\bigcap}
\nc{\bu}{\bigcup}
\nc{\com}{\circ}
\nc{\ssim}{\approx}
\nc{\sq}{\sqrt}
\nc{\rou}{\rho}
\nc{\tf}{\textbf}
\nc{\se}{\simeq}
\nc{\iden}{\mathbb{1}}
\nc{\disu}{\amalg}
\nc{\inclu}{\hookrightarrow}
\nc{\xra}{\xrightarrow}
\nc{\FF}{\mathbb{F}}
\nc{\KK}{\mathbb{K}}
\nc{\NP}{\mathbb{N}_+}
\nc{\emp}{\varnothing}
\nc{\sB}{\mathscr{B}}
\nc{\sR}{\mathscr{R}}
\nc{\Lm}{\Lambda}
\nc{\sm}{\setminus}
\nc{\rra}{\rightrightarrows}
\nc{\op}{\oplus}
\nc{\bop}{\bigoplus}
\nc{\ssum}{\Sigma}

\nc{\diam}{\mm{diam}}
\nc{\diag}{\mm{diag}}
\nc{\msr}{\mathscr}
\nc{\lan}{\langle}
\nc{\ran}{\rangle}
\nc{\rua}{\rightharpoonup}
\nc{\rp}{\RR\mm{P}}
\nc{\PP}{\mathbb{P}}
\nc{\af}{\mathbb{A}}
\usepackage{tikz-cd}
\usetikzlibrary{cd}
\usetikzlibrary{decorations.pathmorphing}
\usepackage{tkz-euclide}
\usetikzlibrary{patterns}
\usetikzlibrary{shapes}
\usetikzlibrary{shapes.misc}
\usetikzlibrary{calc,math}
\usepackage{ifthen}
\usetikzlibrary{calc}
\tikzset{
  my label/.style={font=\scriptsize,inner sep=2pt},
  a/.style={my label,above,node contents={$a$}},
  b/.style={my label,right,node contents={$b$}},
  a-1/.style={my label,above,node contents={$a^{-1}$}},
  b-1/.style={my label,right,node contents={$b^{-1}$}},
}
\newcommand\caley[6]{
  \ifthenelse{0<#1}{
    \pgfmathtruncatemacro\newlev{#1-1}
    \pgfmathtruncatemacro\len{#2}
    \draw[draw=black,-latex] (0,0) -- (\len pt,0) node[pos=.6,#3] coordinate (O);
    \begin{scope}[shift={(O)}]
      \begin{scope}[rotate=90] \caley{\newlev}{\len/2}{#4}{#5}{#6}{#3} \end{scope}
      \begin{scope}[rotate=0]  \caley{\newlev}{\len/2}{#3}{#4}{#5}{#6} \end{scope}
      \begin{scope}[rotate=-90]\caley{\newlev}{\len/2}{#6}{#3}{#4}{#5} \end{scope}
    \end{scope}
  }{\fill[black] circle(1pt);}
}
\usetikzlibrary{lindenmayersystems,arrows.meta}
\newcount\quadrant
\pgfdeclarelindenmayersystem{cayley}{
  \rule{A -> B [ R [A] [+A] [-A] ]}
  \symbol{R}{ \pgflsystemstep=0.5\pgflsystemstep } 
  \symbol{-}{
    \pgfmathsetcount\quadrant{Mod(\quadrant+1,4)}
    \tikzset{rotate=90}
  }
  \symbol{+}{
    \pgfmathsetcount\quadrant{Mod(\quadrant-1,4)}
    \tikzset{rotate=-90}
  }
  \symbol{B}{
    \draw [dot-cayley] (0,0) -- (\pgflsystemstep,0) 
       node [font=\footnotesize, midway, 
         anchor={270-mod(\the\quadrant,2)*90}, inner sep=.5ex] 
           {\ifcase\quadrant$a$\or$b$\or$a^{-1}$\or$b^{-1}$\fi};
    \tikzset{xshift=\pgflsystemstep}
  }
}
\tikzset{
  dot/.tip={Circle[sep=-1.5pt,length=3pt]}, cayley/.tip={Stealth[]dot[]}
}
\nc{\hl}{\colorbox{yellow}}
\nc{\tcr}{\textcolor{red}}
\nc{\sA}{\mathscr{A}}
\nc{\sC}{\mathscr{C}}
\nc{\sD}{\mathscr{D}}
\nc{\sS}{\mathscr{S}}
\nc{\sT}{\mathscr{T}}
\nc{\id}{\mm{id}}
\nc{\Id}{\mm{Id}}
\nc{\opp}{\mm{o}}
\nc{\Ob}{\mm{Ob}}

\nc{\Ad}{\mm{Ad}}
\nc{\uAd}{\underline{\mathrm{Ad}}}
\nc{\wad}{/\uAd}
\nc{\ad}{\mm{ad}}
\nc{\Aut}{\mm{Aut}}
\nc{\Inv}{\mm{Inv}}
\nc{\Gal}{\mm{Gal}}
\nc{\sig}{\sigma}
\nc{\sgn}{\mm{sgn}}
\nc{\Sym}{\mm{Sym}}
\nc{\sym}{\mm{sym}}
\nc{\vphi}{\varphi}
\nc{\vi}{\varphi}
\nc{\vp}{\varpi}
\nc{\vP}{\varPi}
\nc{\bF}{\mathbf{F}}
\nc{\fT}{\mathfrak{T}}
\nc{\fN}{\mathfrak{N}}

\nc{\kap}{\kappa}
\nc{\fX}{\mathfrak{X}}
\nc{\fY}{\mathfrak{Y}}

\nc{\fa}{\mathfrak{a}}
\nc{\fb}{\mathfrak{b}}
\nc{\fB}{\mathfrak{B}}
\nc{\bss}{\mathcal{B}}
\nc{\opn}{\mm{Op}}
\nc{\ops}{\mm{op}}
\nc{\ff}{\mathcal{F}}
\nc{\res}{\tet{res}}
\nc{\OO}{\mathcal{O}}
\nc{\fm}{\mathfrak{m}}
\nc{\fM}{\mathfrak{M}}
\nc{\fn}{\mathfrak{n}}
\nc{\fp}{\mathfrak{p}}
\nc{\fq}{\mathfrak{q}}
\nc{\Mor}{\tet{Mor}}
\nc{\aff}{\mathbb{A}}
\nc{\grass}{\tet{Grass}}
\nc{\sur}{\twoheadrightarrow}

\usepackage{stmaryrd}
\nc{\llb}{\llbracket}
\nc{\rrb}{\rrbracket}
\nc{\cok}{\mm{Coker}}
\nc{\fU}{\mathfrak{U}}
\nc{\wtil}{\widetilde}
\nc{\Ann}{\mathrm{Ann}}
\nc{\fc}{\mathfrak{c}}
\nc{\Ver}{\mm{Ver}}

\nc{\GL}{\mm{GL}}
\nc{\gl}{\mathfrak{gl}}
\nc{\SO}{\mm{SO}}
\nc{\so}{\mathfrak{so}}
\nc{\SL}{\mm{SL}}
\nc{\fsl}{\mathfrak{sl}}
\nc{\mO}{\mm{O}}
\nc{\Bil}{\mm{Bil}}
\nc{\Alt}{\mm{Alt}}

\nc{\gr}{\mm{gr}}
\nc{\End}{\mm{End}}
\nc{\Idem}{\mm{Idem}}
\nc{\lto}{\leadsto}
\nc{\tor}{\mm{Tor}}
\nc{\ext}{\mm{Ext}}
\nc{\lext}{\uli{\mm{Ext}}}
\nc{\loc}{\mm{loc}}
\nc{\Loc}{\mm{Loc}}
\nc{\pr}{\mm{pr}}
\nc{\cO}{\mathcal{O}}
\nc{\cI}{\mathcal{I}}
\nc{\cA}{\mathcal{A}}
\nc{\cD}{\mathcal{D}}
\nc{\cU}{\mathcal{U}}
\nc{\disc}{\mm{disc}}
\nc{\cls}{\mm{cls}}
\nc{\gen}{\mm{gen}}
\nc{\rep}{\mm{rep}}
\nc{\ade}{\mathbb{A}}
\nc{\pslr}{\mm{PSL}_2(\R)}
\nc{\pslz}{\mm{PSL}_2(\Z)}
\nc{\slr}{\mm{SL}_2(\R)}
\nc{\sor}{\mm{SO}_2(\R)}
\nc{\stab}{\mm{Stab}}
\nc{\acts}{\curvearrowright}
\nc{\tang}{\mm{T}^1\mathbb{H}}
\nc{\HH}{\mathbb{H}}
\nc{\x}{\times}
\nc{\xx}{^\times}
\nc{\xs}{^*}
\nc{\xp}{^\perp}
\nc{\Frac}{\mm{Frac}}
\newcommand{\pp}{%
  \mathrel{\ooalign{$\lneq$\cr\raise.22ex\hbox{$\lhd$}\cr}}}
\nc{\ie}{i.e.~}
\nc{\resp}{resp.~}
\nc{\eg}{e.g.~}
\nc{\cf}{cf.~}
\nc{\tiff}{if and only if}
\nc{\enu}{enumerate}
\nc{\hangin}{\hangindent\leftmargini\textup{(i)}~}
\nc{\hangma}{\hspace*{\leftmargini}}
\nc{\indnull}{\mathbb{o}}
\nc{\inv}{^{-1}}
\nc{\du}{^\vee}
\nc{\sups}{\supset}
\nc{\Hom}{\mm{Hom}}
\nc{\iHom}{\mathcal{H}\!\mathit{om}}
\nc{\mat}{\mm{Mat}}
\nc{\Set}{\mm{Set}}
\nc{\Rng}{\mm{Rng}}
\nc{\Grp}{\mm{Grp}}
\nc{\abGrp}{\mm{AbGrp}}
\nc{\LGrp}{\mm{LieGrp}}
\nc{\CLGrp}{\mm{ConnLieGrp}}
\nc{\CSCLGrp}{\mm{CSCLieGrp}}
\nc{\LAlg}{\mm{LieAlg}}
\nc{\TGrp}{\mm{TopGrp}}
\nc{\Diff}{\mm{Diff}}
\nc{\Deck}{\mm{Deck}}
\nc{\Fld}{\mm{Fld}}
\nc{\Top}{\mm{Top}}
\nc{\Ab}{\mm{Ab}}
\nc{\Mod}{\tet{-}\mm{Mod}}
\nc{\Alg}{\tet{-}\mm{Alg}}
\nc{\kFld}{\tet{-}\mm{Fld}}
\nc{\GG}{\mathbb{G}}
\nc{\Sch}{\mm{Sch}}
\nc{\sch}{\mm{sch}}
\nc{\lrs}{\mm{lrs}}
\nc{\red}{\mm{red}}
\nc{\Fun}{\mm{Fun}}
\nc{\sing}{\mm{sing}}
\nc{\dr}{\mm{dR}}
\nc{\G}{\mathbb{G}}
\nc{\bt}{\bullet}
\nc{\bqp}{\oli{\Q}_p}
\nc{\qp}{\Q_p}
\nc{\zp}{\Z_p}

\newcommand{\thrm}[1]{Theorem~\ref{#1}}

\newcommand{\ex}[1]{Example~\ref{#1}}

\newcommand{\conje}[1]{Conjecture~\ref{#1}}
\nc{\dasha}{\dasharrow}
\nc{\Etale}{\'Etale}
\nc{\etale}{\'etale}

\nc{\nil}{\mathfrak{N}}
\nc{\jac}{\mathfrak{J}}
\nc{\rad}{\mm{rad}}
\nc{\lrad}{\mm{l\tet{-}rad}}
\nc{\rrad}{\mm{r\tet{-}rad}}
\nc{\inj}{\rightarrowtail}

\nc{\psupp}{\tet{-supp}}

\nc{\ot}{\otimes}
\nc{\bx}{\wh{\ot}}
\nc{\ha}{\hat{H}}
\nc{\cE}{\mathcal{E}}
\nc{\cF}{\mathcal{F}}
\nc{\cG}{\mathcal{G}}
\nc{\cH}{\mathcal{H}}
\nc{\cL}{\mathcal{L}}
\nc{\cB}{\mathcal{B}}
\nc{\dlim}{\varinjlim}
\nc{\ilim}{\varprojlim}
\nc{\sU}{\mathscr{U}}
\nc{\sV}{\mathscr{V}}
\nc{\sW}{\mathscr{W}}
\nc{\chH}{\check{H}}
\nc{\sI}{\mathcal{I}}
\nc{\psh}{\mm{PSh}}
\nc{\sh}{\mm{Sh}}
\nc{\sep}{\mm{sep}}
\nc{\alg}{\mm{alg}}
\nc{\ab}{\mm{ab}}
\nc{\setmid}{\,|\,}
\nc{\nsetmid}{\!\nmid\!}
\nc{\smid}{\,|\,}
\nc{\bmid}{\,\big|\,}
\nc{\Bmid}{\,\Big|\,}
\nc{\CP}{\mathbb{C}\mm{P}}
\nc{\RP}{\mathbb{R}\mm{P}}
\nc{\HP}{\mathbb{H}\mm{P}}
\nc{\FP}{\mathbb{F}\mm{P}}
\nc{\ltri}{\triangleleft}
\nc{\tri}{\triangle}
\nc{\act}{\curvearrowright}
\nc{\Isom}{\mm{Isom}}
\nc{\Om}{\Omega}
\nc{\om}{\omega}
\nc{\ev}{\mm{ev}}
\nc{\std}{\mm{std}}
\nc{\triv}{\mm{triv}}
\nc{\sob}{\mm{sob}}

\nc{\lip}[1]{\left\| #1 \right\|_\mm{Lip}}
\nc{\cl}{\mm{cl}}
\nc{\Cl}{\mm{Cl}}
\usepackage{upgreek}
\nc{\pres}{\prescript}
\nc{\ttau}{\pres{t}{}\tau}
\nc{\ttl}[1]{\pres{t}{}\tau^{\le#1}}
\nc{\ttg}[1]{\pres{t}{}\tau^{\ge#1}}
\nc{\tl}[1]{\tau^{\le#1}}
\nc{\tg}[1]{\tau^{\ge#1}}
\nc{\trl}{\sT^{\le0}}
\nc{\trll}[1]{\sT^{\le#1}}
\nc{\trg}{\sT^{\ge0}}
\nc{\trgg}[1]{\sT^{\ge#1}}
\nc{\tH}{\pres{t}{}H}
\nc{\tF}{\pres{t}{}F}
\nc{\sr}[1]{#1^\sharp}
\nc{\Fr}{\mm{Fr}}

\makeatletter
\newcommand{\xinclusur}[2][]{%
  \lhook\joinrel
  \ext@arrow 0359\rightarrowfill@ {#1}{#2}%
  \mathrel{\mspace{-15mu}}\rightarrow
}
\makeatother

\makeatletter
\newcommand{\xlinclusur}[2][]{%
  \lhook\joinrel
  \ext@arrow 0359\rightarrowfill@ {#1}{#2}%
  \mathrel{\mspace{-23.25mu}}\leftarrow
}
\makeatother

\makeatletter
\newcommand*{\relrelbarsep}{.386ex}
\newcommand*{\relrelbar}{%
  \mathrel{%
    \mathpalette\@relrelbar\relrelbarsep
  }%
}
\newcommand*{\@relrelbar}[2]{%
  \raise#2\hbox to 0pt{$\m@th#1\relbar$\hss}%
  \lower#2\hbox{$\m@th#1\relbar$}%
}
\providecommand*{\rightrightarrowsfill@}{%
  \arrowfill@\relrelbar\relrelbar\rightrightarrows
}
\providecommand*{\leftleftarrowsfill@}{%
  \arrowfill@\leftleftarrows\relrelbar\relrelbar
}
\providecommand*{\xrightrightarrows}[2][]{%
  \ext@arrow 0359\rightrightarrowsfill@{#1}{#2}%
}
\providecommand*{\xleftleftarrows}[2][]{%
  \ext@arrow 3095\leftleftarrowsfill@{#1}{#2}%
}
\makeatother
\nc{\xrra}{\xrightrightarrows}
\nc{\xlla}{\xrightrightarrows}

\usepackage{enumitem}
\usepackage[colorlinks=true,hyperindex, linkcolor=magenta, pagebackref=false, citecolor=cyan,pdfpagelabels]{hyperref}
\usepackage[capitalize]{cleveref}
\setlist[enumerate]{itemsep=2pt,parsep=2pt,before={\parskip=2pt}}
\usepackage{xurl}

\usepackage{makecell}
\usepackage{multirow}

\usepackage{longtable}
\usepackage{booktabs}
\usepackage{url}
\usepackage{mathrsfs}
\usepackage{boxedminipage}
\usepackage{anyfontsize}
\usepackage{listings}
\usepackage{xcolor}
\newcommand{\RNum}[1]{\uppercase\expandafter{\romannumeral #1\relax}}
\usepackage{bbold}
\usepackage{wrapfig}
\usepackage{multicol}
\usepackage{extarrows}
\usepackage{titletoc}
\usepackage{dynkin-diagrams}
\usepackage{mathtools}
\usepackage{enumitem}
\usepackage[usestackEOL]{stackengine}

\usepackage{pdfpages}
\makeatletter
\DeclareFontFamily{OMX}{MnSymbolE}{}
\DeclareSymbolFont{MnLargeSymbols}{OMX}{MnSymbolE}{m}{n}
\SetSymbolFont{MnLargeSymbols}{bold}{OMX}{MnSymbolE}{b}{n}
\DeclareFontShape{OMX}{MnSymbolE}{m}{n}{
    <-6>  MnSymbolE5
   <6-7>  MnSymbolE6
   <7-8>  MnSymbolE7
   <8-9>  MnSymbolE8
   <9-10> MnSymbolE9
  <10-12> MnSymbolE10
  <12->   MnSymbolE12
}{}
\DeclareFontShape{OMX}{MnSymbolE}{b}{n}{
    <-6>  MnSymbolE-Bold5
   <6-7>  MnSymbolE-Bold6
   <7-8>  MnSymbolE-Bold7
   <8-9>  MnSymbolE-Bold8
   <9-10> MnSymbolE-Bold9
  <10-12> MnSymbolE-Bold10
  <12->   MnSymbolE-Bold12
}{}

\let\llangle\@undefined
\let\rrangle\@undefined
\DeclareMathDelimiter{\llangle}{\mathopen}%
                     {MnLargeSymbols}{'164}{MnLargeSymbols}{'164}
\DeclareMathDelimiter{\rrangle}{\mathclose}%
                     {MnLargeSymbols}{'171}{MnLargeSymbols}{'171}
\makeatother

\numberwithin{equation}{section}
\newtheorem{theorem}{Theorem}
\numberwithin{theorem}{section}
\newtheorem{thm}[theorem]{Theorem}
\newtheorem{prop}[theorem]{Proposition}

\newtheorem{prop/defi}[theorem]{Proposition/Definition}
\newtheorem{lem/defi}[theorem]{Lemma/Definition}
\newtheorem{thm/defi}[theorem]{Theorem/Definition}
\newtheorem{defi/prop}[theorem]{Definition/Proposition}
\theoremstyle{definition}
\newtheorem{defi}[theorem]{Definition}

\newtheorem{rem}[theorem]{Remark}

\newtheorem{exam}[theorem]{Example}
\newtheorem{conj}[theorem]{Conjecture}

\newtheorem{defi/rem}[theorem]{Definition/Remark}

\makeatletter
\patchcmd{\@tocline}
  {\hfil}
  {\leaders\hbox{\,.\,}\hfil}
  {}{}
\makeatother

\makeatletter
    \def\cleardoublepage{\clearpage%
        \if@twoside
            \ifodd\c@page\else
                \vspace*{\fill}
                \hfill
                \begin{center}
                This page is intentionally left blank.
                \end{center}
                \vspace{\fill}
                \thispagestyle{empty}
                \newpage
                \if@twocolumn\hbox{}\newpage\fi
            \fi
        \fi
    }
\makeatother

\usepackage[new]{old-arrows}
\usepackage{verbatim}
\usepackage{bookmark}
\usepackage{scalerel,stackengine}
\stackMath
\newcommand\rwh[1]{%
\savestack{\tmpbox}{\stretchto{%
  \scaleto{%
    \scalerel*[\widthof{\ensuremath{#1}}]{\kern-.6pt\bigwedge\kern-.6pt}%
    {\rule[-\textheight/2]{1ex}{\textheight}}
  }{\textheight}%
}{0.5ex}}%
\stackon[1pt]{#1}{\tmpbox}%
}
\providecommand*{\rightrightarrowsfill@}{%
  \arrowfill@\relrelbar\relrelbar\rightrightarrows
}
\providecommand*{\leftleftarrowsfill@}{%
  \arrowfill@\leftleftarrows\relrelbar\relrelbar
}
\providecommand*{\xrightrightarrows}[2][]{%
  \ext@arrow 0359\rightrightarrowsfill@{#1}{#2}%
}
\providecommand*{\xleftleftarrows}[2][]{%
  \ext@arrow 3095\leftleftarrowsfill@{#1}{#2}%
}
\makeatother

\begin{document}

\title[A new counterexample to Nguyen's conjecture on surface fibration]{A new counterexample to Nguyen's conjecture\\on surface fibration}
\author{Wenyi Cai, Yuanyuan Qian, Hao Xiao, and Lingyu Zhuo}
\address[Previous]{School of Mathematical Sciences, Soochow University, Suzhou, Jiangsu 215006, China}
\email{hxiao@stu.suda.edu.cn}
\address[Current]{Mathematisches Institut, Universit\"at Bonn, Endenicher Allee 60, D-53115, Bonn, Germany}
\email{xiao@uni-bonn.de}
\begin{abstract}
Suppose $f:S\ra\PP^1$ is a surface fibration of genus $g$ with $3$ singular fibers. If two of the singular fibers are semistable, Nguyen conjectured in \cite{MR1648969} that $f$ does not exist for $g\ge2$. However, a counterexample for $g=2$ was discovered by Gong-Lu-Tan in \cite{MR3072164}. Note that such kind of surface fibrations admit strong arithmetic properties but are rare in fact, and as such the counterexamples are important. In this paper, we construct a new one for $g=2$.
\end{abstract}
\maketitle
\tableofcontents

\st{Introduction to surface fibration}

\noindent
This paper mainly concerns surface fibrations, and we review some related notions:
\begin{defi}
Let $S$ be an algebraic surface, and let $C$ be an algebraic curve.
\begin{\enu}[label=\textup{(\roman*)}]\setcounter{enumi}{0}
\item A proper flat morphism $f:S\ra C$ is a {\bf surface fibration} if its fibers are connected curves.
\item If the general fibers of a surface fibration are smooth of genus $g$, then it is of {\bf genus} $g$.
\item A surface fibration of genus $g$ is {\bf elliptic} if $g=1$.
\item A surface fibration $f:S\ra C$ is {\bf isotrivial} if its smooth fibers are all isomorphic.
\item A surface fibration is {\bf smooth} if all fibers are smooth.
\item A surface fibration is {\bf locally trivial} if it is isotrivial and smooth.
\item An locally trivial surface fibration $f:S\ra C$ is {\bf trivial} if there exists an isomorphism $u:S\ra C'\x C$ for some smooth connected curve $C'$ such that $\pr_2\com u=f$.
\end{\enu}
\end{defi}

It is convenient to introduce some numerical invariants of a surface fibration of fixed genus:
\begin{defi}
Let $f:S\ra C$ be a nontrivial surface fibration of genus $g$. Let $\chi(\cO_S)$ be the Euler–Poincar\'e characteristic of the structure sheaf $\cO_S$, $p_g(S)=h^2(\cO_S)=h^0(\om_S)$ the geometric genus of $S$, and $q(S)=h^1(\cO_S)=h^1(\om_S)$ the irregularity of $S$. They satisfy the identity
\[
\chi(\cO_S)=p_g(S)-q(S)+1.
\]
We also have two Chern numbers $c_1^2(S)$ and $c_2(S)$, \ie $c_1^2(S)$ is the self-intersection number of the first Chern class of the tangent bundle and $c_2(S)$ is the degree of the second Chern class of the tangent bundle (which is equal to the topological Euler–Poincar\'e characteristic of $S$). Now we have the {\bf relative numerical invariants} of $f$, \ie
\begin{align*}
K_f^2&:=c_1^2(S)-8(g-1)(g(C)-1),\\
\chi_f&:=\chi(\cO_S)-(g-1)(g(C)-1),\\
e_f&:=c_2(S)-4(g-1)(g(C)-1),\\
\tet{and}\ q_f&:=q(S)-g(C).
\end{align*}
\end{defi}
\begin{rem}
When $f$ is relatively minimal (\ie $S$ has no ($-1$)-curves contained in any fiber of $f$), these invariants are nonnegative integers and satisfy the relative Noether formula $12\chi_f=K_f^2+e_f$.
\end{rem}

\begin{defi}
Let $f:S\ra C$ be a surface fibration, then a singular fiber of $f$ is {\bf semistable} if it is reduced and its only singularities are normal crossings. And a surface fibration is {\bf semistable} if its singular fibers are all semistable.
\end{defi}

There is an important semistable reduction result by Deligne-Mumford:
\begin{thm}
[{\cite[Corollary~2.7]{MR262240}}]
\label{thm:pullBack}
For every surface fibration $f:S\ra C$, there exists a base change $\pi:\wt{C}\ra C$ such that the pullback one $\wt{f}:\wt{S}\ra\wt{C}$ is semistable.
\end{thm}

This indicates that the semistability of a surface fibration is a rather general property. Moreover, in the proof of \cite[Proposition~1.3]{MR321933}, Arakelov showed a striking inequality:
\begin{thm}[Arakelov]\label{thm:Arakelov}
Let $f:S\ra C$ be a nontrivial semistable surface fibration of genus $g\ge2$, and let $s$ be the number of singular fibers of $f$. Then we have
\[
\chi_f\le\fr{g}{2}(2g(C)-2+s).
\]
\end{thm}
\begin{rem}
Note that $s$ is a finite number by \cite[Theorem~12.2.4.(iii)]{MR217086}.
\end{rem}

Subsequently, Vojta proved the canonical class inequality (\cf Liu \cite[Theorem~0.1]{MR1418581}):
\begin{thm}
[{\cite[Corollary~2.2]{MR969124Vojta}}]
\label{thm:canonicalClass}
Let $f:S\ra C$ be a nontrivial semistable surface fibration of genus $g\ge2$, and let $s$ be the number of singular fibers of $f$. Then we have
\[
K_f^2<(2g-2)(2g(C)-2+s).
\]
\end{thm}

Since surface fibrations with small number of singular fibers have some interesting properties, it is worth trying to obtain some classification results. In the first place, there is one for semistable surface fibrations by Beauville:
\begin{thm}
[{\cite[Theorem~2.1.(i)]{MR3618574}}]
\label{thm:singular>=4}
A nontrivial semistable surface fibration $f:S\ra\mathbb{P}^1$ admits at least $4$ singular fibers.
\end{thm}
\begin{rem}
For semistable elliptic fibrations, there is a complete classification list achieved by Beauville in \cite{MR664643}.
\end{rem}
Beauville also made the following conjecture at the end of page 100 in \cite{MR3618574}:
\begin{conj}[Beauville]
\label{conj:Beauville}
A nontrivial semistable surface fibration $f:S\ra\mathbb{P}^1$ of genus $g\ge2$ admits at least $5$ singular fibers.
\end{conj}

Miranda and Persson respectively in \cite{MR1076128} and \cite{MR1069483} finished the classification of semistable elliptic fibrations of rational surfaces over $\PP^1$ with $5$ singular fibers, and there are $10$ classes in total.

Tan proved \conje{conj:Beauville} in \cite{MR1325793} using \thrm{thm:canonicalClass}. Later on, Liu gave another proof in \cite{MR1418581} by investigating the moduli space $\cM_{g,n}$ of Riemann surfaces of genus $g$ with $n$ punctures.

As for the classification of non-semistable surface fibrations, we have not achieved much progress by far due to the lack of powerful tools like \thrm{thm:Arakelov} and \thrm{thm:canonicalClass}.
Given genus $g\ge2$, there are a few results for the classification of surface fibrations over $\mathbb{P}^1$ with given number of singular fibers. We present a classical one derived by Par\v{s}in:
\begin{thm}
[{\cite[Proposition~5]{MR257086}}]
Let $f:S\ra\mathbb{P}^1$ be a nontrivial surface fibration, and let $s$ be its number of singular fibers, then $s\ge2$. If in addition $f$ is non-isotrivial, then $s\ge3$.
\end{thm}

Beauville also constructed a surface fibration with exactly $3$ singular fibers:
\begin{exam}
[{\cite[Proposition~1.1.(ii)]{MR3618574}}]
\label{exam:Beauville}
Let $t\in\PP^1$ and $n\ge3$, then
\[
y^2=x^n-ntx+(n-1)t
\]
is a surface fibration over $\PP^1$ of genus $\lfloor(n-1)/2\rfloor$. It admits exactly $3$ singular fibers at $t=0,1,\infty$. Note that only the singular fiber at $t=1$ is semistable.
\end{exam}

\st{Nguyen's conjecture}

\noindent
It follows from \thrm{thm:singular>=4} that there does not exist a nontrivial semistable surface fibration over $\PP^1$ with only three singular fibers (\cf \ex{exam:Beauville}), and Nguyen made the following conjecture on nontrivial surface fibrations with three singular fibers in the remark following \cite[Theorem~7]{MR1648969}:
\begin{conj}\label{conj:Nguyen}
It seems that there does not exist a nontrivial surface fibration of genus $g\ge2$ over $\mathbb{P}^1$ with three singular fibers and two of them semistable.
\end{conj}
\begin{rem}
If $g=1$, then such kind of surface fibrations do exist, and there are $6$ classes in total classified by Miranda-Persson in \cite[Table~5.2]{MR867347}.
\end{rem}

Nguyen showed that such kind of surface fibrations (if exist) have rather weird but strong properties in the sense of both arithmetic and topology, so he was inclined to negate their existence.

However, this does not hold. The first counterexample was given in \cite[Example~3.1]{MR3072164}.

On the other hand, by Bely\u{\i}'s theorem (see \cite[Theorem~1]{MR534593}), such kind of surface fibrations can be defined on algebraic number fields, and they have interesting applications in number theory. So we are motivated to search new counterexamples.

\st{The new counterexample}

\noindent
Now we present a new counterexample to \conje{conj:Nguyen}:
\begin{prop}
Let $t$ be the coordinate of the base curve $\PP^1$, then
\[
y^2=-(x^4-x^2+t)(x+1)
\]
is a a surface fibration with three singular fibers and two of them semistable.
\end{prop}
\begin{proof}
We denote the surface by $S$.
Note that $S$ is obtained via a double cover of $\mathbb{P}^1\tim\mathbb{P}^1$ which induces a surface fibration $S\ra\mathbb{P}^1$ of genus $2$. Since the discriminant of $-(x^4-x^2+t)(x+1)$ with respect to $x$ is $16t^3(4t-1)^2$, the singular fibers appear at $t=0,1/4,\infty$. From the double cover, we obtain the branch locus of $S\ra\PP^1$ in Figure~\ref{pic:1}.

\begin{figure}[ht!]
\centering
\begin{picture}(400, 150)
\multiput(300, 27)(0, 7.5){16}{\line(0, 1){11}}
\put(291.5, 127.5){\line(1, -1){17}}
\put(331, 120){\makebox(0, 0){$x=-1$}}
\put(318, 60){\oval(36, 36)[l]}
\qbezier(300,45)(302.5,45)(305,47.25)
\put(330, 60){\makebox(0, 0){$x=\infty$}}
\put(303.5, 43){\makebox(0, 0){\tiny $4$}}
\put(300, 10){\makebox(0, 0){$\mathrm{t=\infty}$}}

\multiput(200, 27)(0, 8){16}{\line(0, 1){3}}
\put(230, 140){\makebox(0, 0){$x=0$}}
\put(191.5, 147.5){\line(1, -1){17}}
\put(230, 120){\makebox(0, 0){$x=1$}}
\put(191.5, 127.5){\line(1, -1){17}}
\put(238, 85){\makebox(0, 0){$x=-\fr{\sqrt{2}}{2}$}}
\put(218, 85){\oval(36, 36)[l]}
\qbezier(200,70)(202.5,70)(205,72.25)
\put(203.5, 67.5){\makebox(0, 0){\tiny $2$}}
\put(234, 45){\makebox(0, 0){$x=\fr{\sqrt{2}}{2}$}}
\put(218, 45){\oval(36, 36)[l]}
\put(203.5, 27.5){\makebox(0, 0){\tiny $2$}}
\qbezier(200,30)(202.5,30)(205,32.25)
\put(206, 10){\makebox(0, 0){$\mathrm{t=1/4}$}}

\multiput(100, 27)(0, 8){16}{\line(0, 1){3}}
\put(135, 140){\makebox(0, 0){$x=\infty$}}
\put(91.5, 147.5){\line(1, -1){17}}
\put(136.5, 120){\makebox(0, 0){$x=-1$}}
\put(91.5, 127.5){\line(1, -1){17}}
\put(91.5, 110.5){\line(1, 1){17}}
\put(133, 85){\makebox(0, 0){$x=0$}}
\put(118, 85){\oval(36, 36)[l]}
\qbezier(100,70)(102.5,70)(105,72.25)
\put(103.5, 67.5){\makebox(0, 0){\tiny $2$}}
\put(134, 45){\makebox(0, 0){$x=1$}}
\put(91.5, 52.5){\line(1, -1){17}}
\put(104, 10){\makebox(0, 0){$\mathrm{t=0}$}}
\end{picture}
\vspace*{-4mm}
\caption{The branch locus of the surface fibration $S\ra\PP^1$.}
\label{pic:1}
\vspace*{-2mm}
\end{figure}

According to the Namikawa-Ueno classification of singular fibers of genus $g=2$ in \cite{MR369362}, the singular fibers at $t=0,1/4$ are pictured in Figure~\ref{pic:2}. Similarly, the singular fiber at $t=\infty$ is also pictured in Figure~\ref{pic:2} where $B$ is the rational branch of self-intersection number $-3$ and the others are rational branches of genus $0$ and self-intersection number $-2$.

\begin{figure}[ht!]
\centering
\begin{picture}(300, 120)
\multiput(60, 27)(0, 7.375){12}{\line(0, 1){12}}
\put(70, 92.5){\circle{20}}
\put(70, 52.5){\circle{20}}
\put(60, 10){\makebox(0, 0){$\mathrm{t=1/4}$}}

\multiput(0, 27)(0, 7.25){12}{\line(0, 1){12}}
\put(10, 92.5){\oval(20, 20)[r]}
\put(10, 102.5){\line(-1, 0){20}}
\put(10, 82.5){\line(-1, 0){20}}
\put(10, 52.5){\circle{20}}
\put(0, 10){\makebox(0, 0){$\mathrm{t=0}$}}

\put(120, 65){\line(0, 1){30}}
\put(115, 90){\line(1, 0){50}}
\put(160, 85){\line(0, 1){30}}
\put(115, 110){\line(1, 0){190}}
\put(180, 85){\line(0, 1){30}}
\put(200, 85){\line(0, 1){30}}
\put(195, 90){\line(1, 0){30}}
\put(220, 65){\line(0, 1){30}}
\put(215, 70){\line(1, 0){50}}
\put(260, 45){\line(0, 1){30}}
\put(255, 50){\line(1, 0){50}}
\put(300, 55){\line(0, -1){30}}

\put(117.5, 80){\makebox(0, 0){$\mathrm{1}$}}
\put(140, 94.25){\makebox(0, 0){$\mathrm{2B}$}}
\put(157.375, 100){\makebox(0, 0){$\mathrm{5}$}}
\put(177.125, 100){\makebox(0, 0){$\mathrm{4}$}}
\put(190, 114.25){\makebox(0, 0){$\mathrm{8}$}}
\put(197.125, 100){\makebox(0, 0){$\mathrm{7}$}}
\put(210, 94.25){\makebox(0, 0){$\mathrm{6}$}}
\put(222.75, 80){\makebox(0, 0){$\mathrm{5}$}}
\put(240, 74){\makebox(0, 0){$\mathrm{4}$}}
\put(262.875, 60){\makebox(0, 0){$\mathrm{3}$}}
\put(280, 54.125){\makebox(0, 0){$\mathrm{2}$}}
\put(302.5, 40){\makebox(0, 0){$\mathrm{1}$}}
\put(200, 10){\makebox(0, 0){$\mathrm{t=\infty}$}}
\end{picture}
\vspace*{-4mm}
\caption{The singular fibers of the surface fibration $S\ra\PP^1$.}
\label{pic:2}
\vspace*{-2mm}
\end{figure}

Note that the singular fibers at $t=0,1/4$ are semistable, and the one at $t=\infty$ is of type \RNum{7}*. Therefore, this is indeed a counterexample to \conje{conj:Nguyen}.
\end{proof}

\st*{Acknowledgments}

\noindent
The authors of this paper appreciate the instruction by Prof.~Dr.~Cheng Gong. This paper received financial support from both National University Students' Innovation and Entrepreneurship Training Program, \emph{Research on Security of Digital Signature (201710285014Z)}, and National Natural Science Foundation of China, \emph{Research on Surface Fibration with Three Singular Fibers (11401413)}.

\bookmarksetup{startatroot}
\bibliographystyle{amsalpha}
\bibliography{../algebra}

\end{document}